\documentclass[11pt,a4paper,twoside]{article}

\usepackage[cp1250]{inputenc}
\usepackage{latexsym,amssymb,amsmath}

\title{Simple Moufang loops and Galois extensions}
\date{}
\author{SANDU N. I.}
\begin{document}
\maketitle

\begin{abstract}

The paper establishes an one-to-one correspondence between simple
Moufang loops and Paige loops constructed over Galois extension
over prime field in its algebraic closure. Using this connection
it describes fully  the family of nonassociative finite simple
Moufang loops. It  describes the generators, the loop structure of
subloops, the automorphism group of nonassociative simple Moufang
loops. \vspace*{0.1cm}


\end{abstract}
\vspace*{0.1cm}

Historically, the first  examples of simple Moufang loops were
constructed by Paige in [1]. Namely, given a field $F$, write
$M^{\star}(F)$ for the set of all elements of unit norm in matrix
Cayley-Dickson algebra $C(F)$ over $F$ (also called split octonion
algebra in the literature), and let $M(F)$ be the quotient of
$M^{\star}(F)$ by its centre $Z(M^{\star}(F)) = \{\pm1\}$ (more
detailed definitions will be presented below). Paige proved that
$M(F)$ is a nonassociative simple Moufang loop. These loops are
also known as Paige loops. Liebeck [2] used the classification of
finite simple groups to conclude that there are no other
nonassociative simple Moufang loops besides $M(F)$, $F$ finite.

Further,  [3] (see, also [ 4]) indicates concrete elements $a, b,
c \in M(F)$, that generate the loop $M(F)$, for $F = GF(p)$.
Elements $a, b, c$ contain only $0$ and $1$ in their structure.
This result is generalized in [5], using the classical results on
generators of unimodular groups. It is proved that every Paige
loop $M(F)$, $F = GF(p^n)$, is $3$-generated. The core of the
proof consists in considering a concrete triplet of generators of
$M(F)$. The main result of [6], proved in a quite cumbersome
manner, is that the automorphism group of Paige loop $M(GF(2))$ is
isomorphic to Chevalley group $G_2(GF(2))$. In [7] it is proved
that the automorphism group of the Paige loop $M(F)$, where  $F$
is a perfect field,  is isomorphic to the semidirect product
$G_2(F) \rtimes$ $\textit{Aut}$ $F$, where $G_2(F)$ is the
Chevalley group.

The aforementioned results from [3 -- 7] can be proved by
different methods. But the crucial moment of this proofs relies
heavily on the results of Doro [8], that relate Moufang loops to
groups with triality.  [4] rewrites many proofs using the
geometric loop theory. It also uses the connections between
composition algebras, simple Moufang loops, simple Moufang
$3$-nets, $S$-simple groups and groups with triality. This paper
contains a new approach to investigation of simple Moufang loops -
using the one-to-one connection between simple Moufang loops
(equivalent to Paige loops) and Galois extension over prime field
in its algebraic closure and also normal subgroups of automorphism
group of algebraic closure of prime field. It describes  the
generators, the automorphism group, the subloop structure of
simple Moufang loops. The aforementioned  results regarding the
simple Moufang loops become a particular case of the corresponding
results from this paper.

For basic definitions and properties of Moufang loops, alternative
algebras , fields see [9], [10], [11] respectively. Often the
needed results from this books will be used without reference. It
is also worth mentioning that this paper is a continuation of
paper [12] and, for purposes of completeness, the content of the
latter is included in this paper.

Remind that for an alternative algebra $A$ with  unity the set
$U(A)$ of all invertible elements of $A$ forms a Moufang loop with
respect to multiplication [13].

It is proved by analogy to Lemma 1 from [14, 15].\vspace*{0.1cm}

\textbf{Lemma 1.} \textit{Let $A$ be an alternative algebra with
unity and let $Q$ be a subloop of $U(A)$. Then the restriction of
any homomorphism of algebra $A$  upon $Q$ will be a homomorphism
on the loop. More concretely, any ideal $J$ of $A$ induces a
normal subloop $Q \cap (1 + J)$ of $Q$.} \vspace*{0.1cm}

Let $L$ be a free Moufang loop, let $F$ be a field and let $FL$ be
a loop algebra of loop $L$ over field $F$. We remind that $FL$ is
a free module with basis $\{g \vert g \in L\}$ and the
multiplication of basis elements is defined by their
multiplication in loop $L$. Let $(u,v,w) = uv\cdot w - u\cdot vw$
denote the associator of elements $u, v, w$ of algebra $FL$. We
denote by $I$ the ideal of loop algebra $FL$, generated by the set

$$\{(a,b,c) + (b,a,c), (a,b,c) + (a,c,b) \vert \forall a, b, c \in L\}.$$ It
is shown in [13, 14] that algebra $FL/I$ is alternative and loop
$L$ is embedded (isomor\-phically) in the loop $U(FL/I)$. Further
we identify the loop $L$ with its isomorphic image in $U(FL/I)$.
Hence the free loop $L$ is a subloop of loop $U(FL/I)$. Without
causing any misunderstandings  we will denote by $FL$ the quotient
algebra $FL/I$ and call it alternative loop algebra (in [14, 15]
it is called $''$loop algebra$''$). Sums $\sum_{g\in L}\alpha_gg$,
are elements of algebra $FL$, where $\alpha_g \in F$. Further, we
will identify the field $F$ with subalgebra $F1$ of algebra $FL$,
where $1$ is the unit of loop $L$.

Let now $Q$ be an arbitrary Moufang loop. Then $Q$ has a
representation as a quotient loop $L/H$ of the free Moufang loop
$L$ by the normal subloop $H$. We denote by $\omega H$ the ideal
of alternative loop algebra $FL$, generated by the elements $1 -
h$ ($h \in H$). By Lemma 1, $\omega H$ induces a normal subloop $K
= L \cap (1 + \omega H)$ of loop $L$ and $F(L/K) = FL/\omega H$.

We denote $L/K = \overline{Q}$, thus $FL/\omega H =
F\overline{Q}$. As every element in $FL$ is a finite sum  $\sum_{g
\in L}\alpha_gg$, where $\alpha_g \in F$, $g \in L$, then the
finite sum $\sum_{q \in \overline{Q}}\alpha_qq$, where $\alpha_q
\in F$, $q \in \overline{Q}$ will be elements of algebra
$F\overline Q$.  Let us determine the homomorphism of $F$-algebras
$\varphi: FL \rightarrow F(L/H)$ by the rule
$\varphi(\sum\lambda_qq) = \sum\lambda_qHq$. The mapping $\varphi$
is $F$-linear, then for $h \in H$, $q \in L$ we have $\varphi((1 -
h)q) = Hq - H(hq) = Hq - Hq = 0$. Hence $\omega H \subseteq$
\textit{ker} $\varphi$. The loop $\overline{Q}$ is a subloop of
loop $U(F\overline{Q})$ and as $\omega H \subseteq$ \textit{ker}
$\varphi$, then the homomorphisms $FL \rightarrow FL/\omega H =
F\overline Q$ and $FL \rightarrow FL/$\textit{ker} $\varphi =
F(L/H) = FQ$ induces a homomorphism $\pi$ of  loop $\overline{Q}$
upon loop $Q$. Hence we have. \vspace*{0.1cm}

\textbf{Lemma 2.} \textit{Let $Q$ be an arbitrary Moufang loop.
Then the loop $\overline{Q}$ is embedded in loop of invertible
elements $U(F\overline{Q})$ of alternative algebra $F\overline{Q}$
and the homomorphism $FL \rightarrow FL/\omega H$ of alternative
loop algebra $FL$ induces a homomorphism $\pi: \overline Q
\rightarrow Q$ of loops.}\vspace*{0.1cm}

\textbf{Remark.} The Lemma 2 answers positively to Question 1'
from [16]: is it true that any Moufang loop can be imbedded into a
homomorphic image of a loop of type $U(A)$ for a suitable unital
alternative algebra?\vspace*{0.1cm}

Let now $Q$ be a simple Moufang loop. Then \textit{ker} $\pi$ will
be a proper maximal normal subloop of $\overline Q$. Let $J_1,$
$J_2$ be  proper ideals of algebra $F\overline Q$. We prove that
the sum $J_1 + J_2$ is also a proper ideal of $F\overline Q$.
Indeed, by Lemma 1 $K_1 = \overline Q \cap (1 + J_1),$ $K_2 =
\overline Q \cap (1 + J_2)$ will be  normal subloops of loop
$\overline Q$. We have that $K_1 \subseteq$ \textit{ker} $\pi$,
$K_2 \subseteq$ \textit{ker} $\pi$. Then product $K_1K_2
\subseteq$ \textit{ker} $\pi$, as well. But $K_1K_2 = (\overline Q
\cap (1 + J_1))(\overline Q \cap (1 + J_2)) = \overline Q \cap (1
+ J_1)(1 + J_2) = \overline Q \cap (1 + J_1 + J_2 + J_1J_2) =
\overline Q \cap (1 + J_1 + J_2)$. Hence $\overline Q \cap (1 +
J_1 + J_2) \subseteq$ \textit{ker} $\pi$, i.e. $\overline Q \cap
(1 + J_1 + J_2)$ is a proper normal subloop of $\overline Q$. Then
from Lemma 1 it follows that $J_1 + J_2$ is a proper ideal of
algebra $F\overline Q$, as required.

We denote by $S$ the ideal of algebra $F\overline  Q$, generated
by all proper ideals $J_i$ ($i \in I$) of $F\overline Q$. Let us
show that $S$ is also a proper ideal of algebra $F\overline Q$. If
$I$ is a finite set, then the statement follows from the first
case. Let us now consider the second possible case. The algebra
$F\overline Q$ is generated as a $F$-module by elements $x \in
\overline Q$. Let there be such ideals $J_1, \ldots, J_k$ that for
element $1 \neq a \in \overline Q$ $a \in \sum J_i$ and let us
suppose that for element $b \in \overline Q$ $b \notin \sum J_i$.
We denote by $T$ the set of all ideals of algebra $F\overline Q$,
containing the element $a$, but not containing the element $b$. By
Zorn's Lemma there is a maximal ideal $I_1$ in $T$. We denote by
$I_2$ the ideal of algebra $F\overline  Q$, generated by all
proper ideals of $F\overline  Q$ that don't belong to ideal $I_1$.
Then $S = I_1 + I_2$. $I_1, I_2$ are proper ideals of $F\overline
Q$ and by the first case $S$ is also proper ideal of $F\overline
Q$. By Lemma 1 $K = \overline Q \cap (1 + S)$ is a normal subloop
of $\overline Q$. We denote $\overline{\overline Q} = \overline
Q/K$. Then $F\overline{ \overline Q} = F\overline Q/S$ is a simple
algebra. As $K \subseteq$ \textit{ker} $\pi$ then $\pi$ induce a
homomorphism $\rho: \overline{\overline Q} \rightarrow Q$. Hence
we prove. \vspace*{0.1cm}

\textbf{Lemma 3.} \textit{Let $Q$ be a simple nonassociative
Moufang loop. Then the loop $\overline{\overline Q}$ is embedded
in loop of invertible elements $U(F\overline{\overline Q})$ of
alternative algebra $F\overline{\overline Q} = F(\overline Q)$ and
the homomorphism $FL \rightarrow FL/\omega H$ of alternative loop
algebra $FL$ induces a homomorphism $\rho: \overline{\overline Q}
\rightarrow Q$ of loops.}\vspace*{0.1cm}

Let $F$ be an arbitrary field. Let us consider a classical matrix
Cayley-Dickson algebra $C(F)$. It consists of matrices of form
$\left(
\begin{array}{ll} \alpha_1 & \alpha_{12}\\
\alpha_{21} & \alpha_2  \end{array} \right)$, where $\alpha_1,
\alpha_2 \in F$, $\alpha_{12}, \alpha_{21} \in F^3$. The addition
and multiplication by scalar of elements of algebra $C(F)$ is
represented by ordinary addition and multiplication by scalar of
matrices, and the multiplication of elements of algebra $C(F)$ is
defined by the rule

$$ \left(
\begin{array}{ll} \alpha_1 & \alpha_{12}\\
\alpha_{21} & \alpha_2  \end{array} \right) \left(
\begin{array}{ll} \beta_1 & \beta_{12} \\ \beta_{21} & \beta_2
\end{array} \right) = $$
$$ \left( \begin{array}{ll} \alpha_1 \beta_1 + (\alpha_{12},
\beta_{21}) & \alpha_1 \beta_{12} + \beta_2 \alpha_{12} -
\alpha_{21} \times \beta_{21} \\ \beta_1 \alpha_{21} + \alpha_2
\beta_{21} + \alpha_{12} \times \beta_{12} & \alpha_2 \beta_2 +
(\alpha_{21}, \beta_{12}) \end{array} \right), \eqno{(1)} $$ where
for vectors $\gamma = (\gamma_1, \gamma_2, \gamma_3),$ $\delta =
(\delta_1, \delta_2, \delta_3) \in A^3$  $(\gamma, \delta) =
\gamma_1\delta_1 + \gamma_2\delta_2 + \gamma_3\delta_3$ denotes
their scalar product and $\gamma \times \delta = (\gamma_2\delta_3
- \gamma_3\delta_2, \gamma_3\delta_1 - \gamma_1\delta_3,
\gamma_1\delta_2 - \gamma_2\delta_1)$ denotes the vector product.
Algebra  $C(F)$ is alternative.  For $a \in C(F)$ the norm $n(a)$
is defined by the equality $n(a) = \alpha_1\alpha_2 -
(\alpha_{12}, \alpha_{21})$. If $\alpha_{ij} = (\alpha_{ij}^{(1)},
\alpha_{ij}^{(2)}, \alpha_{ij}^{(3)})$, then
$$n(a) = \alpha_1\alpha_2 - \alpha_{12}^{(1)}\alpha_{21}^{(1)} -
\alpha_{12}^{(2)}\alpha_{21}^{(2)} -
\alpha_{12}^{(3)}\alpha_{21}^{(3)}. \eqno{(2)}$$ The algebra
$C(F)$ is split. Then from [1, Theorem] it follows.
\vspace*{0.1cm}

\textbf{Lemma 4.} \textit{Let $F$ be an arbitrary field. Then the
Moufang loop $M(F) = M^{\star}(F)/<-1>$ of the matrix
Cayley-Dickson algebra $C(F)$ is simple and the loop $<-1>$
coincides with the center of loop $M^{\star}(F)$}. \vspace*{0.1cm}

The loops $M(F)$ of Lemma 4 are sometimes called \textit{Paige
loops}.

By Lemma 4 the center $Z$ of loop $M^{\star}(F)$ coincides with
subloop $<-1>$. As $M^{\star}(F)/Z$ is a simple loop, a question
appears. Is the center $Z$ of loop $M^{\star}(F)$ emphasized by
the direct factor? The answer is negative. Let field $F$ consist
of 5 elements and let $H$ be a direct completion of center $Z$. If
$\alpha$ were the generator of the multiplicative group of field
$F$, then one of the elements $\pm\left(
\begin{array}{ll} \frac{1}{\alpha} & 0 \\ 0 & \alpha \end{array}
\right)$ would lie in $H$. The square of this element is equal to
$-1$, i.e., it lies in the intersection $H \cap Z$, which is
impossible. Therefore center $Z$ cannot have a direct factor in
$M^{\star}(F)$.

Let now $P$ be an algebraically closed field and let $Q$ be a
simple nonassociative Moufang loop. By Lemma 3 the loop
$\overline{\overline Q}$ is embedded in loop of invertible
elements of simple alternative algebra $P\overline{\overline Q}$.
We denote $\overline{\overline Q} = G$.

If $a \in G$, then it follows from the equality $aa^{-1} = 1$ that
$n(a)n(a)^{-1} = 1$, i.e. $n(a) \neq 0$. The associator $(a,b,c)$
of elements $a, b, c$ of an arbitrary loop is defined by the
equality $ab\cdot c = (a\cdot bc)(a,b,c)$. Identity $(xy)^{-1} =
y^{-1}x^{-1}$ holds in Moufang loops. Therefore, if $a, b, c$ are
elements of Moufang loop $G$, then $u = (a,b,c) = (a\cdot
bc)^{-1})(ab\cdot c) = (c^{-1}b^{-1}\cdot a^{-1})(ab\cdot c)$,
$n(u) = n(c^{-1})n(b^{-1})n(a^{-1})n(a)\cdot n(b)n(c) =
n(c)^{-1}n(b)^{-1}n(a)^{-1}n(a)n(b)n(c) = 1$, i.e. $u \in
M^{\star}(P)$. We denote by $G'$ the subloop generated by all
associators of Moufang loop $G$. If $G' = G$, then $G \subseteq
M^{\star}(P)$, i.e. the loop $G$ is embedded in $M^{\star}(P)$.
Now we suppose that $G' \neq G$. It is shown in [17, 18] that the
subloop $G'$ is normal in $G$. The finite sum $\sum_{g \in
G}\alpha_gg$, where $\alpha_g \in P$, $g \in G$ are elements of
algebra $PG$. Let $\eta :PG \rightarrow P(G/G')$ be a homomorphism
of $P$-algebras determined by the rule $\eta(\sum \alpha_gg) =
\sum \alpha_ggG'$ ($g \in G$) and let $P(G/G') = PG/$\textit{ker}
$\eta$. As the quotient loop $P(G/G')$ is non-trivial, then
$PG/$\textit{ker} $\eta \neq PG$. Hence \textit{ker} $\eta$ is a
proper ideal of $PG$. The algebra $PG$ is simple. Then the ideal
\textit{ker} $\eta$ cannot be the proper ideal of $PG$. Hence the
case $G' \neq G$ is impossible and, consequently, the loop $G$ is
embedded in loop $M^{\star}(P)$.

The alternative algebra $PG$ is simple. By Kleinfeld Theorem [10]
it is a Cayley-Dickson algebra over their center. Field $P$ is
algebraically closed. Then algebra $PG$ is split. The matrix
Cayley-Dickson algebra $C(P)$ is also split. But any two split
nonassociative composition algebras over an algebraically closed
field are isomorphic. Therefore algebra $PG = P\overline{\overline
Q}$ is isomorphic to the matrix Cayley-Dickson algebra $C(P)$.
Further we identify $P\overline{\overline Q}$ with $C(P)$. It is
proved. \vspace*{0.1cm}

\textbf{Lemma 5.} \textit{Let $P$ be an algebraically closed field
and let $Q$ be a simple Moufang loop. Then loop
$\overline{\overline{Q}}$ is embedded in loop $M^{\star}(P)$ of
matrix  Cayley-Dickson algebra $P\overline{\overline Q} = C(P)$}.
\vspace*{0.1cm}

It is known that for any field $K$ there exists an algebraic
closure $\overline K$, containing $K$ as subfield. Let $S$ be the
set of all roots of all polynomials of degree $\geq 1$ from
polynomial ring $K[X]$. Then the field $\overline K$ is the adding
of the set $S$ to field $K$.  The elements of $\overline K$ are
polynomials of elements from $S$ with coefficients from $K$.
\textit{Further we will consider that $\Delta$ is a prime field
and denote by $P$ (or $\overline{\Delta}$) its algebraic closure.}
Every field contains as subfield an unique prime field. Then every
subfield of $P$ can be presented  as an extension of prime field
$\Delta$.

Let $\varphi: \overline{\overline Q} \rightarrow M^{\star}(P)$ be
the embedding  considered in Lemma 5. Further we identify
$\overline{\overline Q}$ with $\varphi \overline{\overline Q}$.
The subfield over $\Delta$  of field $P$, generated by matrix
elements $\alpha_i, \alpha_{ij}^{(k)}$, where $\alpha_{ij} =
(a_{ij}^{(1)}, a_{ij}^{(2)}, a_{ij}^{(3)})$,  of all
matrices $\left(\begin{array}{ll} \alpha_1 & \alpha_{12} \\
\alpha_{21} & \alpha_2 \end{array} \right)$ $ \in
\overline{\overline Q}$ will be denoted by $P_Q$ and will be
called \textit{subfield over $\Delta$ of field $P$, corresponding
to matrices loop $\overline{\overline Q}$}. \vspace*{0.1cm}

\textbf{Lemma 6.} \textit{Let $Q$ be a simple Moufang loop and let
$P_Q = F$ be the subfield over $\Delta$ of algebraic closure $P$
corresponding to matrix loop $\overline{\overline Q}$. Then
$M^{\star}(P_Q) = \overline{\overline Q}$.}\vspace*{0.1cm}

\textbf{Proof.} Let $L$ be such a free Moufang loop  that $Q =L/H$
for simple Moufang loop $Q$. Denote by $\widetilde{FL}$ the loop
algebra  and by $FL$ the alternative loop algebra of loop $L$ over
field $F$. Let
$$\widetilde{FL} \rightarrow \widetilde{FL}/I = FL
\rightarrow FL/\omega H = F\overline Q \rightarrow F\overline Q/S
= F\overline{\overline Q}$$ be the homomorphisms of algebras
constructed in Lemmas 2, 3. These homomorphisms of algebras induce
the homomorphisms of loops
$$L \rightarrow L/B = \overline Q \rightarrow \overline Q/K = \overline{\overline Q}.$$

Let now $R= M(F)$ be the simple Moufang loop considered in Lemma 4
and let $L_1$ be such a free Moufang loop  that $R = L_1/H_1$. As
in the case of loop $Q$ we consider  the homomorphisms
$$\widetilde{FL_1} \rightarrow \widetilde{FL_1}/I_1 = FL_1
\rightarrow FL_1/\omega H_1 = F\overline R \rightarrow F\overline
R/S_1 = F\overline{\overline R},$$
$$L_1 \rightarrow L_1/B_1 = \overline R \rightarrow \overline R/K_1 =
\overline{\overline R}.$$ of algebras and loops  respectively.
These homomorphisms induce the homomorphisms $\varphi:
\widetilde{FL_1} \rightarrow F\overline{\overline R}$ and $\psi:
L_1 \rightarrow \overline{\overline R}$ of algebras and loops such
that $\varphi u = \psi u$ for all $u \in L_1$.

By Lemma 5,  $\overline{\overline Q} \subseteq M^{\star}(P_Q) =
\overline{\overline R}$ and $F\overline{\overline Q} = C(F)$ ($F =
P_Q$), where $C(F)$ is matrix Cayley-Dickson algebra. Hence
$F\overline{\overline Q} = F\overline{\overline R}$.  Let us
suppose  that $\overline{\overline Q} \subset \overline{\overline
R}$. Let $a \in \overline{\overline R} \backslash
\overline{\overline Q}$. It  follows from equality
$F\overline{\overline Q} = F\overline{\overline R}$ that $a = \sum
\alpha_ia_i$, where $\alpha_i \in F$, $a_i \in \overline{\overline
Q}$, and the sum $\sum \alpha_ia_i$ is non-trivial, i.e. $i > 1$.
Let $\psi u = a$, $\psi u_i = a_i$ for some $u, u_1 \in L_1$. Then
the inverse image of equality $a = \sum \alpha_ia_i$ under
homomorphism $\varphi: \widetilde{FL} \rightarrow
F\overline{\overline R}$ has a form $u = \sum \alpha_iu_i + \sum
\beta_jv_j$ for some $\beta_j \in F$, $v_j \in L_1$. The sum $\sum
\alpha_i$ in equality $a = \sum \alpha_i$ is non-trivial. Then,
obviously, the sum $\sum \alpha_iu_i + \sum \beta_jv_j$ in
equality $u = \sum \alpha_iu_i + \sum \beta_jv_j$ is also
non-trivial. Hence in loop algebra $\widetilde{FL_1}$ the element
$u \in L_1$  is linearly expressed through the elements of loop
$L_1$. But this  contradicts the definition of loop algebra
$\widetilde{FL_1}$. Consequently, $\overline{\overline Q} =
\overline{\overline R} = M^{\star}(F)$. This completes the proof
of Lemma 6. \vspace*{0.1cm}

Let $Q$ be a simple Moufang loop and let $P_Q$ be the subfield
over $\Delta$ of algebraic closure $P$ corresponding to matrix
loop $\overline{\overline Q}$. From Lemma 6 it follows that
$\overline{\overline Q}$ consists of elements of form $a = \left(
\begin{array}{cc} \alpha_1 & (\alpha_2, \alpha_3, \alpha_4) \\
(\alpha_5, \alpha_6, \alpha_7) & \alpha_8
\end{array} \right)$ with norms one, where $\alpha_1, \ldots,
\alpha_8 \in P_Q$. Denote $T_i =
\{\alpha_i \in P_Q \vert \quad \forall a \in \overline{\overline
Q} \}$, $i = 1, 2, \ldots, 8$. The elements
$$\left(
\begin{array}{cc} \alpha & (\alpha, \alpha, \alpha) \\ (0, 0, 0) & \alpha^{-1}
\end{array} \right), \left(
\begin{array}{cc} \alpha & (0, 0, 0) \\ (\alpha, \alpha, \alpha) & \alpha^{-1}
\end{array} \right), \left(
\begin{array}{cc} \alpha^{-1} & (0, 0, 0) \\ (\alpha, \alpha, \alpha) & \alpha^{-1}
\end{array} \right)$$ are norm one. From here  it follows that $\alpha \in T_i$ implies
$\alpha \in T_j$, $i, j = 1, \ldots , 8$, or $T_1 = T_2 = \ldots =
T_8$. Let us denote $BP_Q = T_i$. Then only the elements from
$BP_Q$ are as components of matrices elements from
$\overline{\overline Q}$. The norms of elements
 $\left(\begin{array}{cc} \alpha & (1, 0, 0) \\ (0, 0, 0) & \alpha^{-1}
\end{array} \right)$,
$\left(
\begin{array}{cc} 1 & (0, \alpha, 1) \\ (0, 0, 0) & 1
\end{array} \right)$ are one for  $\alpha \in BP_Q$.  By (1)
$$\left(
\begin{array}{cc} \alpha & (1, 0, 0) \\ (0, 0, 0) & \alpha^{-1}
\end{array} \right)\left(
\begin{array}{cc} \beta & (1, 0, 0) \\ (0, 0, 0) & \beta^{-1}
\end{array} \right) = \left(
\begin{array}{cc} \alpha \beta & (1, 0, 0) \\ (0, 0, 0) & \alpha^{-1}\beta^{-1}
\end{array} \right),$$
$$\left(
\begin{array}{cc} 1 & (0, \alpha, 1) \\ (0, 0, 0) & 1
\end{array} \right)\left(
\begin{array}{cc} 1 & (0, \beta, 1) \\ (0, 0, 0) & 1
\end{array} \right) =$$
$$= \left(
\begin{array}{cc} 1 + \alpha\beta & (0, \alpha + \beta, 2) \\ (\alpha - \beta, 0, 0) & 1
\end{array} \right).$$

From these equalities it follows that $\alpha, \beta \in BP_Q$
implies $\alpha^{-1}, \alpha\beta, \alpha + \beta, \alpha - \beta
\in BP_Q$. Hence we proved. \vspace*{0.1cm}

\textbf{Lemma 7.} \textit{Let $Q$ be a simple Moufang loop, let
$\Delta$ be a prime field and let $P$ be its algebraic closure.
Then the set $BP_Q$  of $P$  is a subfield over $\Delta$ and $BP_Q
= P_Q$.} \vspace*{0.1cm}

Every field contains an unique simple subfield. Then $\Delta
\subseteq P_Q$ and $P_Q$ is an extension of field $\Delta$,
$\Delta \subset P_Q$.  The field $P$ is algebraically closed. Then
$\Delta \subset P_Q$ is an \textit{algebraically extension}.

We will show that the \textit{extension $\Delta \subset P_Q$ is
normal}. Really, let $\{\alpha, \alpha_1, \break \alpha_2,
\ldots\}$ be a basis of vector space $P_Q$ over $\Delta$. Then
$P_Q = \Delta(\alpha_1, \alpha_2, \ldots)$. Let $f_{\alpha}(X)$ be
an irreducible polynomial of polynomial ring $\Delta[X]$
corresponding to element $\alpha$, $f_{\alpha}(\alpha) = 0$. We
consider that \textit{deg} $f_{\alpha}(X) > 1$. Let $\beta \neq
\alpha$ be a root of $f_{\alpha}(X)$. To  prove the normality of
extension $\Delta \subset P_Q$ it is sufficient to prove that
$\beta \in P_Q$. Let us assume the contrary, that $\beta \notin
P_Q$. We denote $\Delta_{\alpha} = \Delta(\alpha_1, \alpha_2,
\ldots)$. Over field $\Delta_{\alpha}$ the polynomial
$f_{\alpha}(X)$ takes the form $f(X) = (x - \theta_1) \ldots (x -
\theta_r)g_{\alpha}(X)$, where $\theta_1, \ldots , \theta_r \in
\Delta_{\alpha}$ and the polynomial $f_{\alpha}$ is an irreducible
polynomial over $\Delta_{\alpha}$. As $\beta \notin P_Q$ then
\textit{deg} $g_{\alpha} > 1$. If the degree
$|\Delta_{\alpha}(\alpha) : \Delta_{\alpha}| = n > 1$ then by [11,
Prop. 3, cap. VII] and $|\Delta_{\alpha}(\beta) : \Delta_{\alpha}|
= n$. In these cases the elements in $\Delta_{\alpha}(\alpha)$
have the form $\Sigma_{k=0}^{n-1}\theta_k\alpha^k$ and the
elements in $\Delta_{\alpha}(\beta)$ are form
$\Sigma_{k=0}^{n-1}\theta_k\beta^k$, where $\theta_k \in
\Delta_{\alpha}$. Then  the mapping $\varphi$:
$\Delta_{\alpha}(\alpha)\rightarrow \Delta_{\alpha}(\beta):$
$\Sigma_{k=0}^{n-1}\theta_k\alpha^k \rightarrow
\Sigma_{k=0}^{n-1}\theta_k\beta^k$ is an isomorphism over
$\Delta_{\alpha}$.

We have $P_Q = \Delta_{\alpha}(\alpha)$ and $\Delta_{\alpha} \neq
P_Q$. Then $M(\Delta_{\alpha}) \neq M(P_Q) = Q$. From here it
follows that for $\alpha$ there exist such elements $\omega_1,
\ldots , \omega_7 \in \Delta_{\alpha}$ that $\alpha, \omega_1,
\ldots , \omega_7$ are the components of some matrix element $u$
of Moufang loop $M^{\star}(P_Q)$. The norm $n(u) = 1$ or by (2)
$\alpha\omega_1 - \omega_2\omega_3 - \omega_4\omega_5 -
\omega_6\omega_7 = 1$. Then $\varphi(\alpha\omega_1 -
\omega_2\omega_3 - \omega_4\omega_5 - \omega_6\omega_7) =
\varphi(1)$, $\varphi(\alpha)\omega_1 - \omega_2\omega_3 -
\omega_4\omega_5 - \omega_6\omega_7 = 1$, $\beta\omega_1 =
\omega_2\omega_3 + \omega_4\omega_5 + \omega_6\omega_7 + 1$,
$\beta\omega_1 \in P_Q$. By Lemma 7 $\beta \in P_Q$. We get a
contradiction as $\beta \notin P_Q$. Hence $\beta \in P_Q$.
Consequently, any root of polynomial $f_{\alpha}(X)$ lies in
$P_Q$. Then $f_{\alpha}(X)$ factors over $\Delta$ in  $P_Q$ as a
product of linear factors. Consequently, $P_Q$ is a splitting
field of the set of irreducible polynomials $\{f_{\alpha}(X) \vert
\alpha  = \alpha_1, \alpha_2, \ldots\}$ over $\Delta$ . By [11,
Theorem 4, cap. VII] the extension $\Delta \subset P_Q$ is normal.

The prime field $\Delta$ can be only the  field of rational
numbers $R$ or the Galois field $GF(p)$. In the first  case
\textit{char} P = 0 and by [19] the extension $\Delta \subset P_Q$
is perfect. Recall, a field $F$ is called \textit{perfect} if any
irreducible polynomial of polynomial ring $F[X]$ is separable. Let
us consider the second  case. By [19] the field $GF(p)$ is perfect
and any algebraical extension of perfect field is perfect. Hence
the field $P_Q$ over $\Delta$ is perfect. In particular, the
extension $\Delta \subset P_Q$ is separable. The algebraically
separable and normal extension is called \textit{Galois}.
Consequently, the extension $\Delta \subset P_Q$ is Galois. Hence
we proved. \vspace*{0.1cm}

\textbf{Lemma 8.} \textit{The extension $\Delta \subset P_Q$ is
Galois and $P_Q$ is a perfect field.} \vspace*{0.1cm}

\textbf{Theorem 1.} \textit{Let $\Delta$ be a prime field and let
$P$ be its algebraic closure. Only and only the Paige loops
$M(F)$, where $F$ is a Galois extension over $\Delta$ in $P$, are
with precise till isomorphism nonassociative simple Moufang
loops.}\vspace*{0.1cm}

\textbf{Proof.} Let $Q$ be a simple Moufang loop. By Lemma 7 the
subfield $P_Q = BP_Q$ over $\Delta$ corresponds to it. By Lemma 8
the extension $\Delta \subset P_Q$ is Galois. If $Q$, $G$ are
simple Moufang loops and $Q \neq G$, then $P_Q = BP_Q \neq BP_G =
P_G$.

Conversely, let $F$ be a Galois extension over $\Delta$ of $P$ and
let $P_Q$ be the subfield of $P$ corresponding to  Paige loop $Q =
M(F)$. By Lemma 8 $P_Q$ is a normal extension over $\Delta$. It is
clearly that $P_Q \subseteq F$. We suppose that $1 \neq a \in F
\backslash P_Q$. Let $f(X)$ be an irreducible polynomial of
polynomial ring $\Delta[X]$ of element $a$, $f(a) = 0$. Let $a,
a_1, \ldots, a_n$ be the roots of $f(X)$ and let $S$ be the
splitting field of polynomial $f(X)$ over $\Delta$. By [11,
Theorem 3, cap. VII] the normality of some extension $K$ over
$\Delta$ equals to the fact, that every irreducible polynomial of
$\Delta[X]$, with a root in $K$, factorizes into linear factors in
$K$. The extensions $F, P_Q$ over $\Delta$ are normal. Then $\{a,
a_1, \ldots, a_n\} \subseteq F$, $\{a, a_1, \ldots, a_n\}
\nsubseteq P_Q$ $M(S) \subseteq M(F) - Q$ and $\{a_1, a_2, \ldots
, a_n\} \nsubseteq P_Q$,  $M(S) \nsubseteq M(P_Q) = Q$. We get a
contradiction. Hence $F = P_Q$. This completes the proof of
Theorem 1. \vspace*{0.1cm}

\textbf{Theorem 2.} \textit{Let $\Delta$ be a prime field, let $P$
be its algebraic closure and let $q = p^n$, where $p = p(Q)$ is a
prime, $n = n(Q)$ is an integer. Then for a nonassociative Moufang
loop $Q$ the following statements are equivalent:}

\textit{1) $Q$ is a finite simple loop;}

\textit{2) $Q$ is isomorphic to one of the Paige loop $M(F)$,
where $F \subseteq P$ is a Galois field $GF(q)$ over $\Delta$;}

\textit{3) $Q$ is isomorphic to Paige loop $M(F)$, where $F
\subseteq P$ is a splitting field  over $\Delta$ of  $(q - 1)$-th
roots of unity;}

\textit{4) $Q$ is isomorphic to Paige loop $M(F)$, where $F$ is
the finite field $GF(q)$.}\vspace*{0.1cm}

\textbf{Proof.} It is clear that the loop $Q$ is finite then and
only then the subfield $BP_Q$ over $\Delta$ from Lemma 7 is
finite. Hence $BP_Q = GF(q)$, i.e. 1) $\Rightarrow$ 2). The
implication 2) $\Rightarrow$ 1) follows from Theorem 1.

$2) \Leftrightarrow 3)$. The Galois field $GF(q)$ is defined in a
unique manner in algebraic closure  $P$ as a splitting field of
polynomial $x^q - x$ and its elements are the roots of this
polynomial. But the set of all roots of polynomial $x^q - x$
consists from $\{0\}$ and the set of roots of polynomial $x^{q -1}
- 1$. Hence $2) \Leftrightarrow 3)$.

Finally, the prime field $\Delta$ is fixed in the conditions of
the theorem, but arbitrary. Let \textit{char} $\Delta = p$ and let
$GF(p^n)$ be a Galois field over $\Delta$. Then $\Delta \subseteq
GF(p^n)$. In such a case the equivalence of 1) and 4) derives from
the equivalence of 1) and 2). This completes the proof of Theorem
2. \vspace*{0.1cm}

Let $\Delta$ be a prime field and let $P$ be its algebraic
closure. By Theorem 1 only and only the Paige loops $M(F)$, where
$F$ is a Galois extension over $\Delta$ in $P$, are with precise
till isomorphism nonassociative simple Moufang loops. By Theorem 2
with precise till isomorphism  the only nonassociative finite
simple Moufanf loops are the loops $M(F)$, where $F = GF(p^n)$ is
a Galois field over $\Delta$ in $P$. Next let us consider the
simple Moufanf loops in form of $M(F)$, $F \subseteq P$. It is
clear that for such loops the inclusion relation $\subseteq$ is
meaningful. \vspace*{0.1cm}

\textbf{Corollary 1.} \textit{1. Only the minimal nonassociative
finite simple Moufang loops are the Paige loops $M(GF(p))$, $p$
prime. Only the minimal nonassociative infinite simple Moufang
loops are the Paige loops $M(Q)$, where $Q$ is the field of
rational numbers.}

\textit{2. Only the  simple Moufang loops $M(\Delta)$, where
$\Delta$ is a prime field, do not contain proper nonassociative
subloops.}

\textit{3. The loop $M(GF(q^m))$ is contained into loop
$M(GF(p^n))$ when and only when $q = p$ and $n$ is  divisible by
$m$.}

\textit{4. Only the  maximal nonassociative simple Moufang loops
is the Paige loops $M(P)$, where $P$ is an algebraic closure of
some prime field $\Delta$. For various prime fields $\Delta$ the
corresponding loops $M(P)$ are isomorphic.}\vspace*{0.1cm}

\textbf{The proof} follows immediately   from Theorems 1, 2 and
well known fact from the field theory:  $GF(p^n) \subseteq
GF(p^m)$ when and only when $m$ is divisible by $n$.
\vspace*{0.1cm}

\textbf{Corollary 2}. \textit{If three elements $a, b, c$ of
simple Moufang loop $M(F)$ are not connected through associative
law, for example, $ab\cdot c \neq a \cdot bc$, and the matrix
components of these elements  $a, b, c$ generate the field $F$, in
particular, at least one component is a primitive of field $F$,
then these elements $a, b, c$ generate the loop $M(F)$.}
\vspace*{0.1cm}

Any element $\neq 0$ is a primitive element in prime field
$GF(p)$. Then from Corollary 1 it follows. \vspace*{0.1cm}

\textbf{Corollary 3}. \textit{Let $\Delta = GF(p)$. Then any three
elements $a, b, c \in M(\Delta)$, not connected through the
associative law, generate the loop $M(\Delta)$.} \vspace*{0.1cm}

Now using the Theorem 1 and Lemma 8 we generalize the main result
of [7], stated at the beginning of this article.\vspace*{0.1cm}

\textbf{Proposition 1.}  \textit{The automorphism group of
arbitrary Paige loop $M(F)$ is isomorphic to the semidirect
product $G_2(F) \rtimes G(F/\Delta) = G_2(F) \rtimes$
$\textit{Aut}$ $F$,  where $G_2(F)$ is the Chevalley group,
 $\textit{Aut}$ $F = G(F/\Delta)$ is the Galois
group of Galois extension of a field $F$ over a prime field
$\Delta$ in algebraic closure $\overline{\Delta}$.}\vspace*{0.1cm}

\textbf{Proof.} It is known that the automorphism groups  of prime
fields $\Delta$ coincide with the identical automorphism,
$\textit{Aut}$ $\Delta = \{id\}$. Then $G(F/\Delta) =
\textit{Aut}$ $F$. \vspace*{0.1cm}

From Proposition 1 it follows.\vspace*{0.1cm}

\textbf{Corollary 4.} \textit{The automorphism group of Paige
group $M(F)$ is isomorphic to the Chevalley group $G_2(F)$ if and
only if the field $F$ is prime.}\vspace*{0.1cm}

The automorphism group of Galois field $GF(p^n)$ is the cyclic
group of order $n$ and is generated by automorphism $\varphi: a
\rightarrow a^p$. If by item 3) of Corollary 1 $GF(p^m) \subseteq
GF(p^{mk}$ then the set of all automorphisms of \textit{Aut}
$GF(p^{mk})$, which induces the identical automorphism on
$GF(p^m)$ forms a cyclic group generated by automorphism
$\varphi^k$. Then from Proposition 1 it follows. \vspace*{0.1cm}

\textbf{Corollary 5.} \textit{The automorphism group of  Paige
loop $M(GF(p^n))$ is isomorphic to the semidirect product $G_2(F)
\rtimes$ $\textit{Aut}$ $GF(p^n)$. The group \break $\textit{Aut}$
$GF(p^n)$ is described above.}\vspace*{0.1cm}

Let $K$ be a field,  let $G =$ \textit{Aut} $K$ and let $H$ be a
subgroup of group $G$. \textit{The fixed field of group $G$} is
the subfield $K^H$ of $K$ defined by $K^H = \{x \in K \vert
\varphi x = x \quad \forall \varphi \in H\}$. \textit{The Galois
group} $G(K/E)$ of normal extension $K$ over $E$ is the set of all
field automorphisms of  $K$ which keep all  elements of $E$ fixed.
Let $F$ be a subfield, $E \subset F \subset K$. In [11, Theorem 3,
cap. VIII] it is proved that the extension $F$ over $E$ is normal
when and only when  the subgroup $G(K/F)$ is normal in $G$. In
such a case $G(F/E), \cong G/G(K/F)$. This isomorphism induces an
one-to-one mapping between the set of the normal extension $S$
over $E$ in $K$ and the set of all normal subgroups $H = G(K/F)$
of $G$, get by formula $S = K^H$.

We consider now a normal extension $\Delta \subset K$, where
$\Delta$ is a prime field. If $\Delta$ is the field of rational
numbers then the extension $\Delta \subset K$ is separable [11,
19] and hence $\Delta \subset K$ is a Galois extension. Let now
$\Delta$ be the Galois field $GF(p)$. We have \textit{Aut} $K =
G(K/\Delta) = G$ and $K^G = \Delta$ as \textit{Aut} $F =
\{\textit{id}\}$ if and only if $F$ is a prime field. Then from
[11, Proposition 11, cap. V; 19] it follows that the normal
extension $\Delta \subset K$ is separable. Hence and in case
$\Delta = GF(p)$ the normal extension $\Delta \subset K$ is
Galois.

The Theorem 1 identifies all nonassociative simple Moufang loops
with Paige loops $M(F)$ constructed over all Galois extensions $F$
over prime field $\Delta$ in it algebraic closure
$\overline{\Delta}$. By Lemma 7 the field $F$ determines the loop
$M(F)$ in a unique manner. \vspace*{0.1cm}

\textbf{Theorem 3.} \textit{Let $\Delta$ be a prime field, let
$\overline{\Delta}$ be its algebraic closure and let $P$ be a
Galois extension of $\Delta$ in $\overline{\Delta}$. We denote by
$\Sigma$  the set of all Galois extensions of $\Delta$ in $P$.
Then between the set $\Sigma_L$ of all Paige loops $M(F)$, where
$F \in \Sigma$, and the set $\Sigma_G$ of all normal subgroups $H
= G(P/F)$ of automorphism group $G =$ \textit{Aut} $P$ there
exists an one-to-one mapping defined by $F = P^H$.  Concretely:}

\textit{1)  each Paige loop $M(F), F \in \Sigma$, has a
corresponding  normal subgroup $H = G(P/F) \in \Sigma_G$ induced
by isomorphism \textit{Aut} $F \cong$ \textit{Aut} $P/G(P/H)$;}

\textit{2) the field $F$ is defined in a unique manner $F = P^H$;}

\textit{3) for each normal subgroup $H \in \Sigma_G$ we can find a
field $F \in \Sigma$ that has the above described relationship
with the normal subgroup $H$;}

\textit{4) let $M(F), H$ and $M(F^{\prime}), H^{\prime}$ be the
pairs described in item 1. Then $M(F) \subset M(F')$ (resp. $M(F)
\neq M(F')$ or $M(F) \supset M(F')$) if and only if $H \supset H'$
(resp. $H \neq H'$ or $H \subset H'$).}\vspace*{0.1cm}

Let $\Delta$ be a prime field and let $\overline{\Delta}$ be their
algebraic closure. We denote by $\mathcal M(\overline{\Delta})$
the set of all subloops of Paige loop $M(\overline{\Delta})$,
which are  nonassociative simple loops (i.e. Paige loops $M(F)$,
$F \subset \overline{\Delta}$), by $\mathcal R(\overline{\Delta})$
we denote the set of all Galois extensions of subfield $F$ of
$\overline{\Delta}$ over $\Delta$ and by $\mathcal
G(\overline{\Delta})$ we denote the set of all normal subgroups of
automorphism group \textit{Aut} $\overline{\Delta}$. The sets
$\mathcal M(\overline{\Delta})$, $\mathcal R(\overline{\Delta})$,
$\mathcal G(\overline{\Delta})$ form respectively the lattices
$\boldmath L(\mathcal M(\overline{\Delta}))$, $\boldmath
L(\mathcal R(\overline{\Delta}))$, $\boldmath L(\mathcal
G(\overline{\Delta}))$ with respect to intersection $A \cap B$ as
intersection of sets $A, B$ and union $A \cup B$ as the least
algebra that contains  subalgebras $A, B$. This latices are full
with zero $M(\Delta)$, $\Delta$, unitary group and with unity
$M(\overline{\Delta})$, $\overline{\Delta}$, \textit{Aut}
$\overline{\Delta}$ respectively (more detail see Theorem 2).

Let $\Delta_1$, $\Delta_2$ be a prime fields and let $Q$ be a
simple Moufang loop. By Theorem 1 $Q \cong M(F_1)$ and $Q \cong
M(F_2)$ for some Galois extensions $\Delta_1 \subset F_1 \subset
\overline{\Delta_1}$ and $\Delta_2 \subset F_2 \subset
\overline{\Delta_2}$. Further, from Theorem 3 it follows that the
intersection of two normal subgroups of group \textit{Aut}
$\overline{\Delta}_1$ corresponds to  the union of fields,
corresponding to these subgroups and union of normal subgroups
corresponds to the intersection of fields. Then from here and
definition of Paige loop it follows.\vspace*{0.1cm}

\textbf{Proposition 2.} \textit{Let $\Delta_1$, $\Delta_2$ be a
prime fields. Then the lattices $\boldmath L(\mathcal
M(\overline{\Delta_1}))$, $\boldmath L(\mathcal
M(\overline{\Delta_2}))$, $\boldmath L(\mathcal
R(\overline{\Delta_1}))$, $\boldmath L(\mathcal
R(\overline{\Delta_2}))$ are isomorphic among themselves and are
inverse isomorphic with isomorphically lattices $\boldmath
L(\mathcal G(\overline{\Delta_1}))$, $\boldmath L(\mathcal
G(\overline{\Delta_2}))$.}\vspace*{0.1cm}

It is known [20] that the lattice of all normal subgroups of an
arbitrary group is Dedekind's (modular), i.e. satisfies the
identity $x \cap (y \cup z) = y \cup (x \cap z)$. The dedekind
property is maintained  under inverse isomorphisms. Then from
Proposition 1 it follows.\vspace*{0.1cm}

\textbf{Corollary 6.} \textit{The lattices considered in
Proposition 2 are Dedekind's.}\vspace*{0.1cm}

Let $Q, G$ be  isomorphic loops, let $\sigma: Q \rightarrow G$ be
an isomorphism an let $\sigma^{-1}: G \rightarrow Q$ be the
inverse isomorphism. Denote by \textit{Iso} $(Q, G)$ the set of
all isomorphisms from $Q$ on $G$. If $\varphi \in Aut$ $Q$ then
$\sigma\varphi \in Iso$ $(Q, G)$. Analogically, if $\psi \in Aut$
$B$ then $\psi\sigma \in Iso$ $(G, Q)$. If $\tau \in Iso$ $(Q, G)$
is another isomorphism then $\tau^{-1}\sigma \in Aut$ $Q$,
$\tau\sigma^{-1} \in Aut$ $G$. Hence two isomorphisms  differ to
automorphism. Moreover, the automorphism groups $Aut$ $Q$, $Aut$
$G$ are isomorphic with respect to one-to-one mappings $\varphi
\rightarrow \sigma\varphi\sigma^{-1}$, $\psi \rightarrow
\sigma^{-1}\psi\sigma$. Consequently, the group $Aut$ $Q$ is
defined uniquely with an accuracy to mapping, analogically to
conjugate.

According  to Theorem 1 and to the aforementioned we will consider
that any simple Moufang loop $Q$ is a Paige loop $M(F)$, where $F$
is a Galois extension over a prime field $\Delta$ in algebraic
closure $\overline{\Delta}$, for described the group $Aut$ $Q$. In
this form  $M(F)$ the automorphism group \textit{Aut} $Q$ is
described in Proposition 1. For the finite loop $Q$ the group
\textit{Aut} $Q$ is described in Corollary 5 more detailed.

Now we pass to infinite loop $Q$. Many properties of finite Galois
groups are proved using the calculation  methods, but for infinite
Galois groups it is used, as usually the topological methods (see,
for example, [21]). To recall this. Let $E \subset K$ be a Galois
extension. Íŕ Galois group $G(K/E) = G$ is defined the Krull
topology, taking as the fundamental system of neighborhoods of
unity the set of subgroups $G(S/E)$, where $S$ ranges over all
intermediate fields, $E \subset S \subset K$, such that  $E
\subset S$ is a  finite Galois extension. We prove that the family
of normal subgroups $G(E/k)$, such that  the extension $k \subset
E$ is finite, defines the same topology. From here it follows that
the group $G$ can be presented in the form of  \textit{inverse
limit}
$$G = \lim_{\leftarrow}G(S_i/E),$$ where $S_i$ ranges over
all intermediate fields such that $E \subset S_i$ is a finite
Galois extension. Hence $G$ is a profinite group. We remind that
that a topological group is called \textit{profinite} if is
isomorphic to the inverse limit of an inverse system of discrete
finite groups. The class of all profinite groups and the class of
all compact and totally disconnected topological  groups coincide.
From here it follows that the group $G$ is compact and totally
disconnected.

Finite groups are profinite in  discrete topology and for a prime
field $\Delta$  $G(F/\Delta) =$ \textit{Aut} $F$. Then from the
aforementioned, Propositions 2, 3 it follows.\vspace*{0.1cm}

\textbf{Theorem 4.} \textit{Let $Q$ be a simple Moufang loop. Then
$\textit{Aut}$ $Q \cong G_2(F) \rtimes$ $\textit{Aut}$ $F$, where
$G_2(F)$ is the Chevalley group, $\textit{Aut}$ $F = G(F/\Delta)$,
$G(F/\Delta)$  is the Galois group of Galois extension of a field
$F$ over a prime field $\Delta$ in algebraic closure of $\Delta$.
The group \textit{Aut} $F$ is profinite and if the field $F$ is
infinite then the group $\textit{Aut}$ $F$ is compact and totally
disconnected in Crull topology, and is isomorphic to the  inverse
limit of an inverse system of finite normal subgroups of group
$\textit{Aut}$ $F$.}
\smallskip

\bigskip

\begin{flushleft}
Sandu Nicolae Ion,\\[3mm]

Tiraspol State University of Moldova,\\

Iablochkin str. 5,\\

Kishinev MD-2045, Moldova\\[3mm]

E-mail: sandumn@yahoo.com \\
\end{flushleft}
\end{document}